\documentclass{amsart}
\usepackage{amssymb,psfrag,url,textcmds}
\usepackage[dvips]{graphicx}

\newcommand{\R}{\ensuremath{\mathbb{R}}}
\renewcommand{\d}{\partial}

\theoremstyle{plain}

\newtheorem{theorem}{Theorem}[section]
\newtheorem{corollary}[theorem]{Corollary}
\newtheorem{lemma}[theorem]{Lemma}
\newtheorem{proposition}[theorem]{Proposition}

\theoremstyle{definition}

\newcommand{\sfrac}[2]{{#1/#2}}

\newcommand{\diam}[1]{\textrm{diam}({#1})}

\begin{document}
\title[ ]{When Soap Bubbles Collide}
\author[ ]
{Colin Adams, Frank Morgan and 
John M. Sullivan}

\maketitle

\section{Introduction.} \label{S:Introduction}

\begin{figure}[htbp]
\begin{center}
\includegraphics*[width=2in]{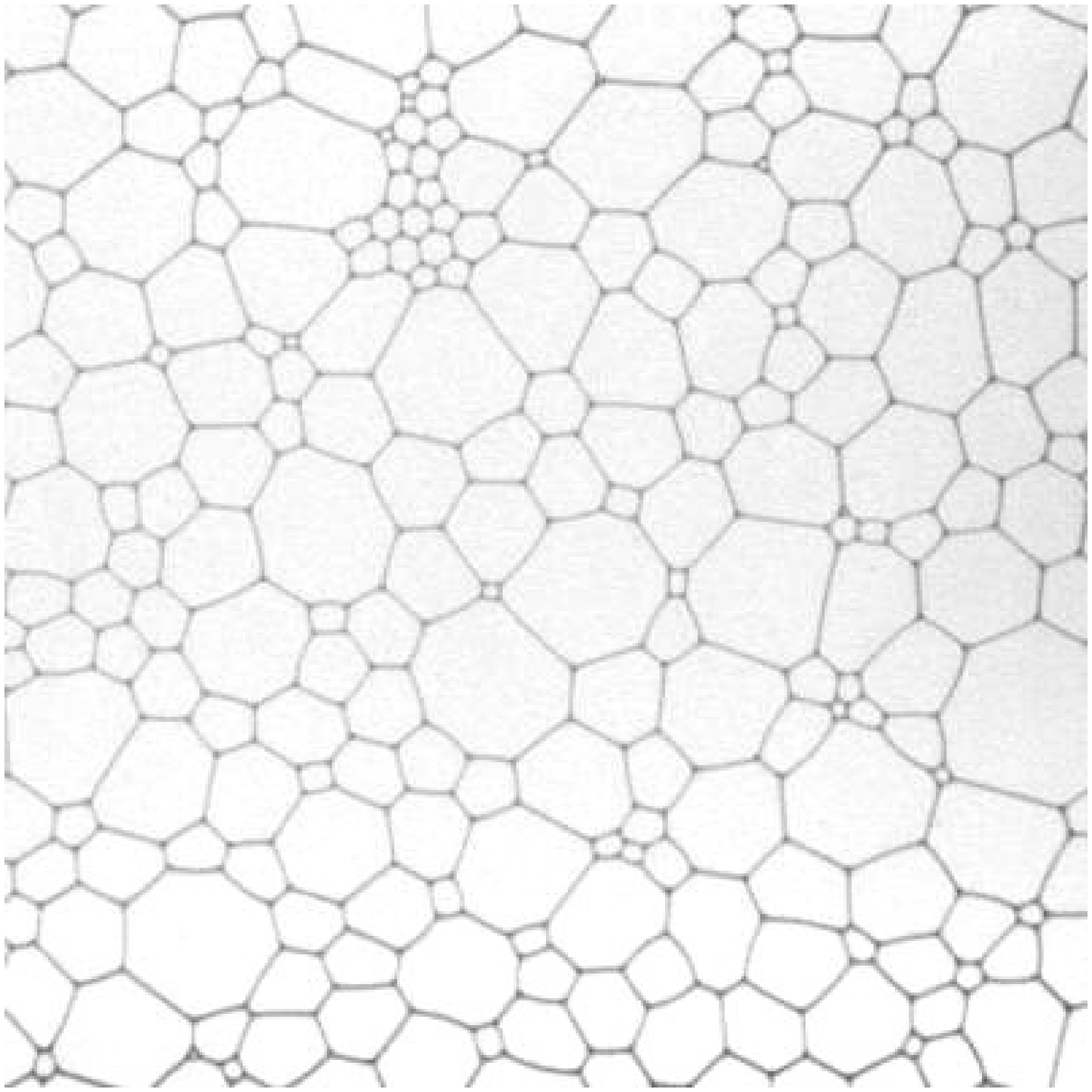}
\hspace{0.5in}
\includegraphics*[width=2in]{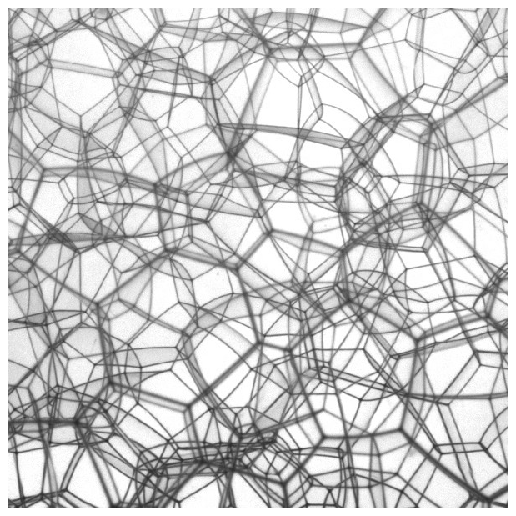}
\end{center}
\caption{\label{Fig:foams} Planar soap bubbles meet in
threes (photograph, left, taken by Olivier Lorderau at Rennes), while
decompositions of $\R^3$ like the soap froth on the right
typically meet in fours (photograph, right, taken
by Sigurdur Thoroddsen at the University of Illinois).}
\end{figure}

\begin{figure}[htbp]
\begin{center}
\includegraphics*[width=1.8in]{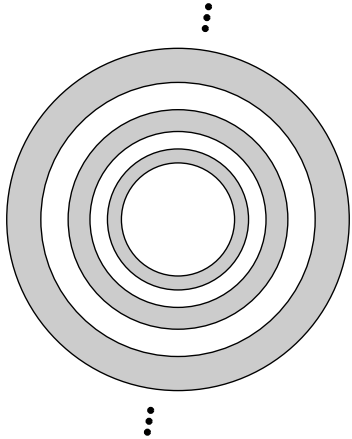}
\end{center}
\caption{\label{Fig:annulic} Nested circles decompose the plane into
pieces of finite volume (a disk and infinitely many annuli)
that meet only in twos.  Revolving this figure around a line
through the center would give a similar example in three dimensions.
This decomposition can be two-colored, as shown.}
\end{figure}

Planar soap bubble froths, as in Figure~\ref{Fig:foams}, have bubbles
meeting at most three at a time. Of course, we can create other
locally finite decompositions of the plane
into closed subsets with disjoint interiors in which arbitrarily many
of these pieces meet at a point. But are
there decompositions that meet at most two at a point? The answer is yes, as
for example occurs with a disk and concentric annuli as in
Figure~\ref{Fig:annulic}. For sets of bounded
diameter, however, the answer is no, as we will see shortly.

\begin{figure}
\includegraphics*[width=1.8in]{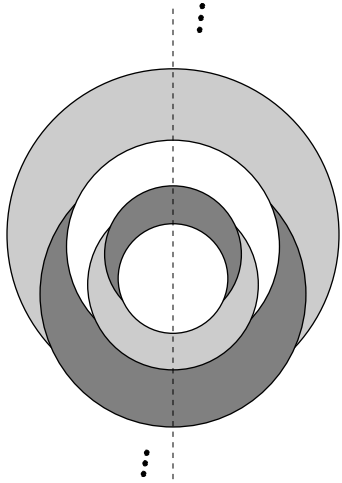}
\caption{\label{Fig:pancakes} The ``pancakes'' in this decomposition (of the plane, or,
 when revolved around the dashed line, of space) meet only in threes.  In
fact, it can be globally three-colored, as shown.}
\end{figure}

Similarly, pieces in a decomposition of $\R^3$ typically meet
in fours, as in the right-hand panel of Figure~\ref{Fig:foams}. Again, we can 
create examples meeting only
in twos, such as a ball and concentric spherical shells. Even when we
restrict ourselves to decompositions into pieces that are topologically balls,
there are examples whose pieces meet at most in threes.
For instance, start with a ball. Cover all but a south polar cap with a
pancake. Now cover all but the pancake's north polar cap with a second,
southern pancake. Continue layering pancakes, alternating over the top and
under the bottom, to fill $\R^3$. See Figure~\ref{Fig:pancakes} for the
two-dimensional version, which can be rotated about a north-south axis
to obtain the version in three (or higher) dimensions.

If we restrict ourselves to pieces of bounded diameter,
however, this cannot occur.
In fact, there are no decompositions of any $\R^n$
into pieces of bounded diameter meeting at most $n$ at a time.
This can be seen from the following version of Lebesgue's covering
theorem (compare [9, Thm. IV 2]):

\begin{theorem}\label{T:Mainresult}
The unit $n$-cube $C^n=[0,1]^n$ has no finite decomposition into pieces
of diameter less than~$1$ that meet at most $n$ at a point.
\end{theorem}

We will sketch a proof using the Brouwer fixed point theorem
in section~\ref{S:Proofs}.  Note that this result is sharp, in the sense that
the unit $n$-cube does admit a piecewise smooth decomposition into
sets of diameter at most $1+\epsilon$ that meet at most~$n$ at a point.
For example, for $n=3$, take the product of an interval with
a fine two-dimensional hexagonal honeycomb in the unit square.
This yields pieces meeting at most three at a time.
Similarly for general~$n$, take the product of an interval
with a fine generic decomposition of the $(n-1)$--cube.
     
We consider locally finite decompositions (where any compact subset meets only finitely many pieces of the decomposition).

\begin{corollary}\label{C:manifold}
For any Riemannian $n$-manifold $M$ there is an~$s$ such that
any locally finite decomposition into pieces of diameter at most~$s$ must have
a point where at least $n+1$ pieces intersect.
\end{corollary}

\begin{proof} A homeomorphic copy of the Euclidean unit $n$-cube $C^n$
can always be mapped into~$M$ with bounded distortion of distance.
Any decomposition of $M$ into pieces with sufficiently small diameters
will yield a decomposition of $C^n$ into pieces of 
diameter less than~$1$. Hence the pieces of the original decomposition of 
$M$ must, in places, intersect $n+1$ at a time.
\end{proof} 
 
Any decomposition of~$\R^n$ into pieces of uniformly bounded diameter
could be rescaled to give a decomposition into pieces of diameter
less than~$s$, so we get as a corollary the result mentioned earlier:

\begin{corollary}\label{C:Rn}
When $n \ge 1$, there is no locally finite decomposition of $\R^n$ into sets of bounded
diameter meeting at most $n$ at a point. 
\end{corollary}

Similarly, by inscribing an $n$-cube
inside the $n$-ball $B^n=\{x \in \R^n: ||x|| \le 1\}$, we also obtain the following result:

\begin{corollary}\label{C:balls}
There is no decomposition of the
unit $n$-ball $B^n$ into sets of diameter less than $\sfrac{2}{\sqrt{n}}$ that
meet at most $n$ at a point. \qed
\end{corollary}

\begin{figure}
\includegraphics*[width=2.2in]{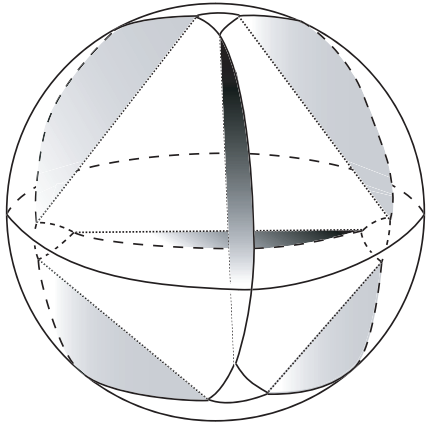}
\caption{\label{Fig:simplex4} A decomposition of $B^3$ into five
regions meeting at most three at a time.
To make the decomposition, we use the four faces of a regular
tetrahedron---slightly larger than the one inscribed in the ball---together
with six ``wings'' outwards from the edges to the sphere.}
\end{figure}

\noindent  We conjecture that this constant $\sfrac{2}{\sqrt{n}}$ is not sharp.
The best subdivision for the unit $n$-ball of which we are aware is obtained 
as follows in the case $n=3$.  Inscribe a regular tetrahedron $\sigma_3$ in the unit $3$-ball, and then (by scaling)
pull its vertices just outside the ball, shaving off that portion of the
simplex now outside the ball.  What remains of the simplex will be one of the five pieces in our decomposition.  Its diameter is just over
$\sqrt{8/3}$.
Then divide the exterior of the simplex in the ball into four pieces,
one for each face of the simplex, by radially projecting the edges of the simplex
out to the sphere, sweeping out four planar surfaces that decompose the remainder of the ball into four pieces, each of diameter just under $\sqrt{2+2/ \sqrt{3} }$. See Figure~\ref{Fig:simplex4}.
The five pieces meet at most three at a point.

Similarly, we can decompose the unit $n$-ball into $n+2$ pieces, one of
diameter just over $\sqrt{2 + \sfrac{2}{n }}$ and the rest of
diameter just under $\sqrt{2+\sqrt{2-2/n}}$.
 The $n+2$ pieces meet at most $n$ at a point.
 
In the case $n=2$, the four pieces in this decomposition have equal
diameters $\sqrt3$, and in fact we prove
in Proposition~\ref{P:equilateral} that this is the minimal-diameter 
decomposition into pieces meeting at most two at a time.
When $n>2$, the truncated simplex can be expanded until its diameter matches 
the diameters of the other pieces. This appears to 
yield the minimal-diameter decomposition into pieces
meeting at most~$n$ at a time.

\subsection*{Implications for soap bubble clusters.}
Planar soap bubble clusters as in
Figure~\ref{Fig:foams} are known to meet in at most threes [10], and 
triple junctions are always at angles of $120$ degrees.
Similarly, soap bubble clusters in three dimensions
meet in at most fours [10]. Corollary~\ref{C:balls} implies that if the
soap bubble cluster covers a ball of diameter at least $\sfrac{\sqrt{3}}{2}$
times the largest diameter of a bubble in the cluster, then there must be
a point where four bubbles meet.

Soap bubble clusters in four dimensions often meet in fives,
but another kind of singularity, where bubbles meet in eights---as
in the cone over the two-skeleton of the four-cube---is also allowed
for soap films [3].
Hence the results here  give only lower bounds on the number of bubbles
of bounded diameter that must meet in higher dimensions.

Unfortunately, there are no known useful \emph{a priori} bounds on the
diameters of bubbles in a cluster.  Even for a fixed number~$m$
of unit volumes in~$\R^n$, there may be many area-minimizing clusters
(all with the same global minimum for total surface area).
At least we can give the following bound on cluster
diameter. Actually, similar arguments show that the space of such
clusters modulo translations is smoothly compact. No such compactness,
or even existence, is known for infinite clusters or for variable,
unbounded numbers of bubbles.

\begin{proposition} For fixed $m$ and $n$,
area-minimizing clusters of $m$ unit volumes
in $\R^n$ have uniformly bounded diameter.
\end{proposition}

\begin{proof}[Sketch of proof]
Take any sequence of minimizing clusters.
As in [10, Theorem 13.4], using translation we may assume that they
converge weakly. (By regularity, a minimizing cluster must be
connected, so that no repetitions are necessary to recover all the
volume.) By a limit argument (see [2] or [11, proof of Theorem 2.1]) the
mean curvatures are weakly bounded.
By monotonicity ([2, Section 5.1] or [10, Theorem 9.3 and subsequent remark])
the diameters are bounded.
\end{proof}

\begin{figure}
\includegraphics*[width=2.0in]{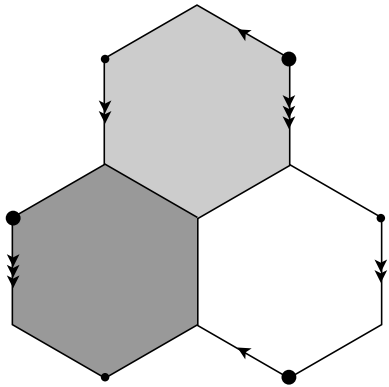}
\caption{\label{Fig:hex} Minimizing double bubbles in flat two-tori
include hexagonal tilings [6], but this can never happen
in $n$-tori for $n > 2$.  Here we see a hexagonal torus divided into
three congruent regions.  (The torus is obtained from the figure by
identifying the boundary as marked; the three vertices marked with
small dots get identified, as do the three with large dots.)}
\end{figure}

The Double Bubble Problem in a Riemannian manifold seeks the
least-area way to enclose and separate two regions of prescribed volumes.
If the manifold is closed, then the complement of the double bubble is a third region of prescribed volume. According to a recent theorem of Corneli et al.~[6], in flat
two-tori there are five types of solutions, including the hexagonal
tiling of Figure~\ref{Fig:hex}.
By a \emph{tiling}, we mean a cluster that, when lifted to the
universal cover, gives a partition into regions of finite volume
or, equivalently, finite diameter.
Such tilings posed the main difficulty in [6].
In the three-torus, Cari{\'o}n et al.~[5] conjectured that
there are no such minimizing double-bubble tilings.
That conjecture (and more) follows easily from Lebesgue's theorem:

\begin{corollary}\label{C:notilings}
In a compact Riemannian $n$-manifold whose universal cover is~$\R^n$,
tilings cannot occur as minimizing $m$-bubble clusters when $m < n$.
\end{corollary}

\begin{proof}
In the universal cover, such a tiling would decompose a large ball into
sets of bounded diameter meeting at most $m+1\le n$ at a time.
\end{proof}

More generally, the Lyusternik--Shnirel{\tprime}man category (see [8])
gives a lower bound on the number of contractible components of a
decomposition of a manifold. The proof of Corollary~\ref{C:notilings}
implies that if the universal cover
of a compact $n$-manifold contains large Euclidean balls, then the
Lyusternik--Shnirel{\tprime}man category is $n+1$ (the largest possible value).

\subsection*{Dimension theory.}
Initially, the concept of the dimension
of a space was not carefully delineated and was inferred from the number of
coordinates required to describe a point. However, at the end of
the nineteenth century, Cantor proved that there was a one-to-one
correspondence between the points in a plane and the points on a line.
Moreover, Peano proved there was a continuous map of an interval onto
a square. It then became clear that an appropriate formal definition of dimension
was necessary. In 1913, L.E.J. Brouwer introduced a definition of
dimension for a topological space
that has evolved into the definition used today.
Topologists characterize $n$-dimensional spaces
by the property that open covers can be refined to meet at most
$n+1$ at all points (see [9]).  Lebesgue's Theorem 1.1 (actually due
to Brouwer [4]) shows that the topological dimension of $\R^n$
is at least~$n$.

\section{Proofs.} \label{S:Proofs}

We now outline the proof of Theorem~\ref{T:Mainresult},
following [9, Section IV.3], which the reader can consult for details.
First we need a lemma equivalent to the Brouwer fixed-point theorem.

\begin{lemma}\label{L:cube}
Let $(F_i,F'_i)$ be the $n$ pairs of opposite faces of the unit $n$-cube $C^n$.
If $B_1,\dots,B_n$ are subsets of the cube such that~$B_i$
separates~$F_i$ from~$F'_i$ in $C^n$, then their
intersection $\bigcap B_i$ is nonempty.
\end{lemma}

\begin{proof}  By definition,  $C^n \smallsetminus B_i$ is the 
disjoint union of two relatively open sets $U_i$ and~$U'_i$
containing $F_i$ and~$F'_i$ respectively.  For each $x \in C^n$
define a vector $v(x)$ whose $i\text{th}$ component
is given by $\pm d(x, B_i)$, where $d(x, B_i)$ signifies
the distance from $x$ to $B_i$
and we use $+$ if $x$ belongs to $U_i$ and $-$ if $x$ lies in $U'_i$.
Define the continuous map~$f$ by $f(x) := x+v(x)$.
This map takes $C^n$ into itself, since the $i \text{th}$ coordinate
changes by the distance to~$B_i$, which is less that the distance
to the opposite face ($F_i$ or $F'_i$).

By the Brouwer fixed point theorem there must be a point~$x_0$
fixed by~$f$. In other words, $v(x_0)=0$ is the zero vector, so
$x_0$ lies in the intersection $\bigcap B_i$.
\end{proof}

The second, easier lemma, which we state without proof, will allow
us to build the sets $B_i$ needed for the first lemma. 

\begin{lemma}\label{L:extend}
Suppose $X\subset Y$ is closed, and $F,F'\subset Y$
are disjoint and closed.  If $A\subset X$ is a closed
set separating $F\cap X$ from $F'\cap X$ in~$X$, then there
exists a closed set $B\subset Y$ separating $F$ from~$F'$
such that $B\cap X\subset A$.
\qed\end{lemma}

\begin{proof}[Sketch of proof of Theorem~\ref{T:Mainresult}]
Suppose there were a decomposition of the unit cube~$C^n$ into pieces of
diameter less than~$1$ meeting at most $n$ at a time, and let~$F_i$ and~$F'_i$
be the faces of the cube, as in Lemma 2.1. Let~$L_1$ be the union of those 
decomposition pieces intersecting $F_1$.
Next, let~$L_2$ be the union of those pieces that intersect~$F_2$ but that
were not already used to build~$L_1$.
Continue in this manner to define $L_3,\ldots,L_n$.
Finally, let~$L_{n+1}$ be the union of the remaining pieces,
those that touch none of $F_1,\ldots,F_n$.

Now take~$K_1$ to be the boundary of~$L_1$ as a subspace of the cube
(equivalently, $K_1 := L_1 \cap \big(\bigcup_{i=2}^{n+1}L_i\big)$).
Then $K_1$ separates~$F_1$ from~$F'_1$ in~$C^n$.
Next, within~$K_1$ as ambient space, let $K_2$ be the boundary
of $L_2\cap K_1$;
it separates $K_1 \cap F_2$ from $K_1 \cap F'_2$ in $K_1$.
Define $K_3, \dots, K_n$ similarly, noting that $K_n=\bigcap_{i=1}^{n+1} L_i$. 
By Lemma~\ref{L:extend}, we can extend the sets $K_i\subset K_{i-1}$ to sets
$B_i\subset C^n$ separating~$F_i$ from~$F'_i$; the intersection of
the~$B_i$ is contained in~$K_n$.
By Lemma~\ref{L:cube}, this intersection (and hence~$K_n$) is nonempty.
Since no piece was used in the construction of more than one~$L_i$,
this proves that $n+1$ of the pieces from the decomposition intersect.
\end{proof}

We next want to derive, in the case $n=2$, the sharp estimate
for the diameter of pieces in a partition of the unit $n$-ball meeting
at most $n$ at a point.  We first quote a standard result from convex
geometry, a combination of the separation theorem
(which is one version of the Farkas alternative)
and Carath\'eodory's theorem.
(See for instance [12, Theorems 1.3.4 and 1.1.4].)

\begin{lemma}\label{L:vectors}
For a subset $X$ of $\R^n$, the following conditions are equivalent:
\begin{itemize}
\item $X$ lies in no open halfspace bounded by a plane through the origin;
\item the origin is in the convex hull of~$X$;
\item the origin is in the convex hull of some subset $X'$ of $X$
  containing at most $n+1$ points.
\end{itemize}
\qed
\end{lemma}

We prove our next proposition for arbitrary dimensions,
though we will use it later only for $n=2$.

\begin{proposition}\label{prop:sphere}
Among all subsets of the unit $(n-1)$--sphere
that are not contained in any open hemisphere,
the set of vertices of the regular $n$-simplex uniquely minimizes diameter.
\end{proposition}

\begin{proof}
Let $A_n$ be the set of $n+1$ vertices of the regular $n$-simplex
inscribed in the unit sphere~$S^{n-1}$.  Every pair of points
in~$A_n$ realizes the diameter, which (in the 
metric inherited from~$\R^n$) is $d_n := \sqrt{2+2/n}$.
It is easy to see (by induction) that $A_n$ is the only
set of more than $n$ points in which all pairs have equal distance.

The proposition is of course trivial for $n=1$, since any proper
subset of $A_1=S^0$ is contained in a hemisphere.

We now proceed by induction.  Suppose we have $B\subset S^{n-1}$
contained in no open hemisphere and with diameter at most $d_n$.
By Lemma~\ref{L:vectors}, we can find a finite subset~$C\subset B$
of $N\le n+1$ points, still contained in no open hemisphere.

Suppose $C$ is contained in the closed hemisphere bounded
by~$S$, some great~$S^{n-2}$.  Then $C\cap S$ is contained
in no open hemisphere of~$S$, for otherwise we could
tilt~$S$ to get an open hemisphere of the original $S^{n-1}$
which contains all of~$C$.  Thus by induction
$$\diam{B}\ge \diam{C}\ge \diam{C\cap S} \ge d_{n-1} > d_n,$$
contradicting the choice of~$B$.
Noting that any subset of at most $n$ points in~$S^{n-1}$
does lie in some closed hemisphere, we conclude that $N=n+1$.

We claim that $C=A_n$.  From this, it follows that $B=A_n$.
For suppose $B$ included an additional point $p$;
since $p$ lies within (the radial projection to the sphere of)
some face of the simplex,
its distance to the opposite vertex would be greater than $d_n$.

To prove the claim, consider the family of all configurations
of $n+1$ (not necessarily distinct) points on~$S^{n-1}$ that are
not contained in any open hemisphere.  Since this family
is compact, diameter can be minimized, and since $d_{n-1}>d_n$,
minimizing configurations have $n+1$ distinct points.
We may assume that $C$ is such a minimizer, and in fact
is one in which the diameter is realized by as few pairs as possible.
We note that $\diam{C}>\sqrt2$, since otherwise $C$ would be contained
in the closed hemisphere around any of its points.

If $C\ne A_n$, then there is some point $p_0\in C$
which is at distance $\diam{C}$ from some but not
all of the other points.  Let $p_1, \ldots, p_k$
be those at maximal distance, with $k\le n-1$.
To check that we can move $p_0$ slightly to~$p'_0$
such that the distance to each $p_i$ decreases, consider
the stereographic projection $S^{n-1}\smallsetminus\{p_0\}\to\R^{n-1}$.

The points $p'_0\ne p_0$ for which the distance to $p_i$ has decreased
project to an open halfspace, which includes the
origin (the image of $-p_0$) since the original
distance was more than $\sqrt2$.  The intersection
of at most $n-1$ such halfspaces is necessarily
unbounded, meaning there are choices for~$p'_0$
arbitrarily near to~$p_0$.  Thus we can decrease
all $k$ distances while maintaining the hemisphere condition.

The new configuration either has lower diameter or has fewer pairs
realizing its unchanged diameter; in either case this contradicts
the choice of the minimizer~$C$.
\end{proof}

\begin{proposition}\label{P:equilateral}
An inscribed equilateral triangle (Figure 6) provides a
least-diameter smooth decomposition of the open unit disc
into relatively closed sets that meet at most two at a point.
\end{proposition}
 
\begin{proof}
Suppose that there is a smooth decomposition into closed
pieces that meet at most two at a point, with smaller diameters than the
inscribed equilateral triangle decomposition.

\begin{figure}
\includegraphics*[width=1.8in]{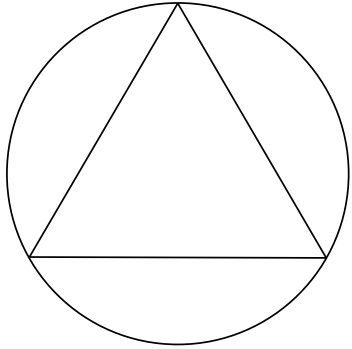}
\caption{\label{Fig:equilateral}An inscribed equilateral triangle provides a least-diameter decomposition of the open unit disk into regions meeting in twos.}
\end{figure}

We may assume that each piece is connected; if not we simply replace each
disconnected piece by the collection of its components.
We may further assume that the closure of each piece touches the unit
circle.   For otherwise, an innermost piece~$A$ would be a topological disk
with a single boundary loop~$\d A$.
Along that loop, there are no triple points,
meaning that $A$ is completely surrounded by a single piece~$B$.  We can then
simply absorb $A$ into~$B$.  This does not increase the diameter, since
if $p\in A$ realizes the maximum distance to any point of $A\cup B$,
then $p\in \d A\subset B$.
Finally, we may assume that the boundary of each piece is
made up alternately of arcs of the circle and chords of the circle,
as replacing arbitrary curves by chords increases no diameters.

An outermost chord cannot cut off a piece that
has boundary intersecting the circle in a
length greater than $\pi$, for the diameter of this piece would be too
large. By absorbing outermost regions into adjacent regions of larger
diameters, we can assume that the closure of precisely one piece~$S$
has disconnected intersection with the circle. We may assume that
the rest are chordal sections. By Proposition~\ref{prop:sphere},
$S$ has diameter at
least as great as the inscribed equilateral triangle, a contradiction.
\end{proof}

\subsection*{ACKNOWLEDGEMENTS}
We wish to thank Joseph Corneli and Stephen Hyde
for their inspiring questions, Branko Gr\"{u}nbaum for
referring us to Lebesgue's theorem, and the editor
for detailed and helpful comments on the submitted version.


\medskip

\noindent {\bf COLIN ADAMS} received his B.S. from MIT and his Ph.D. from
the University of Wisconsin in 1983. He has spent the vast majority of his
career at Williams College. His research interests include knot theory
and hyperbolic 3-manifold theory.  Recipient of the MAA Haimo Teaching
Award in 1998 and the Robert Foster Cherry Great Teachers Award in 2003,
he is the author of  the textbook ``The Knot Book'' and the comic book
``Why Knot?'', and co-author of the humorous supplements to calculus
``How to Ace Calculus'' and ``How to Ace the Rest of Calculus''. He is
also the humor columnist for the Mathematical Intelligencer.

\noindent {\it Bronfman Science Center, Williams College, Williamstown, MA 01267 USA}

\noindent {\it Colin.Adams@williams.edu}

\medskip

\noindent {\bf FRANK MORGAN} and his students love soap bubbles,
minimal surfaces, and the calculus of variations. Founding director
(with Adams) of the Williams College NSF ``SMALL'' Undergraduate Research
Project, inaugural recipient of the MAA Haimo Teaching Award and past
vice-president of the MAA, he has published over a hundred articles and
six books. His MAA Math Chat book grew out of his Math Chat call-in
TV show and column at MAA online. His most recent books are texts on
``Real Analysis'' and ``Real Analysis and Applications.''

\noindent {\it Bronfman Science Center, Williams College, Williamstown, MA 01267 USA}

\noindent {\it Frank.Morgan@williams.edu}

\medskip

\noindent {\bf JOHN SULLIVAN}  received his Ph.D. from Princeton in 1990,
following earlier degrees from Harvard and Cambridge.
After teaching six years each at the Universities of
Minnesota and Illinois, he moved in 2003 to the Technical
University of Berlin.  His research concerns geometric
optimization problems, including bubbles and other surfaces
of constant mean curvature, as well as geometric knot theory
and discrete differential geometry.  His mathematical artwork
(computer-generated prints, sculptures and videos) has been
exhibited in Bologna, Boston, London, New York and Paris, among other places.

\noindent {\it TU Berlin, MA 3--2, Str.~des 17.~Juni~136,
D--10623 Berlin, Germany}

\noindent{\it Sullivan@Math.TU--Berlin.DE}

\end{document}